\documentclass[12pt,a4paper,reqno]{amsart}
\usepackage{amsmath}
\usepackage{amsfonts}
\usepackage{amssymb}
\usepackage{bm} 

\newcommand\eps{\varepsilon}
\newcommand\R{{\mathbf{R}}}

\newcommand\Z{{\mathbf{Z}}}

\renewcommand\Re{{\operatorname{Re}}}

\def \endprf{\hfill  {\vrule height6pt width6pt depth0pt}\medskip}

\theoremstyle{plain}
  \newtheorem{theorem}[subsection]{Theorem}
  
  \newtheorem{proposition}[subsection]{Proposition}
  \newtheorem{lemma}[subsection]{Lemma}
  \newtheorem{corollary}[subsection]{Corollary}

\theoremstyle{remark}

\theoremstyle{definition}

\include{psfig}

\begin{document}

\title[Counterexample to bilinear endpoint Strichartz]{A counterexample to an endpoint bilinear Strichartz inequality}
\author{Terence Tao}
\address{Department of Mathematics, UCLA, Los Angeles CA 90095-1555}
\email{tao@math.ucla.edu}

\vspace{-0.3in}
\begin{abstract}
The endpoint Strichartz estimate $\| e^{it\Delta} f \|_{L^2_t L^\infty_x(\R \times \R^2)} \lesssim \|f\|_{L^2_x(\R^2)}$ is known to be false by the work of Montgomery-Smith \cite{montgomery-smith}, despite being only ``logarithmically far'' from being true in some sense.  In this short note we show that (in sharp constrast to the $L^p_{t,x}$ Strichartz estimates) the situation is not improved by passing to a bilinear setting; more precisely, if $P, P'$ are non-trivial smooth Fourier cutoff multipliers then we show that
the bilinear estimate
$$\| (e^{it\Delta} P f) (e^{it\Delta} P' g)  \|_{L^2_t L^\infty_x(\R \times \R^2)} \lesssim \|f\|_{L^2_x(\R^2)}
\|g\|_{L^2_x(\R^2)}
$$
fails even when $P$, $P'$ have widely separated supports.
\end{abstract}

\maketitle

\section{Introduction}

Consider the Schr\"odinger propagation operators $e^{it\Delta}$ in two spatial dimensions, defined on 
$L^2_x(\R^2)$ via the Fourier transform as
$$ \widehat{e^{it\Delta} f}(\xi) := e^{-4\pi^2 it|\xi|^2} \hat f(\xi)$$
where of course $\hat f$ is the Fourier transform of $f$,
$$ \hat f(\xi) := \int_{\R^2} e^{-2\pi i x\cdot \xi} f(x)\ dx.$$
As is well known, the function $u(t,x) := e^{it\Delta} f(x)$ is the unique (distributional) solution in $C^0_t L^2_x(\R \times \R^2)$ to the free Schr\"odinger equation $iu_t + \Delta u = 0$ with initial datum $u(0,x) = f(x)$.  The problem of understanding the spacetime distribution of $u$ in terms of the $L^2_x$ norm of $f$ is thus of importance in the theory of both linear and nonlinear Schr\"odinger equations.  One fundamental family of estimates in this direction are
the \emph{Strichartz estimates}, which in this two-dimensional setting assert that
\begin{equation}\label{etf}
 \| e^{it\Delta} f \|_{L^q_t L^r_x(\R \times \R^2)} \lesssim_{q,r} \|f\|_{L^2_x(\R^2)}
\hbox{ whenever } \frac{1}{q}+\frac{1}{r} = \frac{1}{2} \hbox{ and } 2 < q \leq \infty,
\end{equation}
where we use $X \lesssim Y$ to denote the estimate $X \leq CY$ for some constant $C$, and more generally use $X \lesssim_{q,r} Y$ to denote $X \leq C_{q,r} Y$ for some constant $C_{q,r}$ depending only on $q,r$, and so forth.
The first estimate of this type (in the case $q=r$) appeared in \cite{strich}, and has since led to many generalisations and extensions; see for instance \cite{tao:keel} and the references therein for further discussion
and for a (short) proof of \eqref{etf}.

Note that the endpoint case $(q,r) = (2,\infty)$ just barely fails to verify the hypotheses of \eqref{etf}.  In \cite{montgomery-smith} this estimate was in fact shown to be false, even with frequency localisation:

\begin{theorem}[Linear endpoint Strichartz fails]\label{counter}\cite{montgomery-smith} Let $P$ be a Fourier multiplier with symbol in $C^\infty_0(\R^2)$ (thus $\widehat{Pf} = \varphi \hat f$ for some $\varphi \in C^\infty_0(\R^2)$)
which is not identically zero.  Then there does not exist a constant $C > 0$ for which one has the estimate
$$ \| e^{it\Delta} P f \|_{L^2_t L^\infty_x(\R \times \R^2)} \leq C\|f\|_{L^2_x(\R^2)}$$
for all $f \in L^2_x(\R^2)$.
\end{theorem}

The proof of this theorem proceeds via a counterexample constructed via Brownian motion; for the convenience of the reader we briefly sketch a proof of this result in an appendix.

While the endpoint Strichartz estimate fails, there are many senses in which the estimate is ``almost'' true.  For instance, in \cite{tao:radial}, \cite{stefanov} it was shown that the estimate becomes true again if one restricts $f$ to be spherically symmetric, or more generally if one performs a suitable averaging in the angular variable.  Also, in dimensions three and higher $d \geq 3$ the endpoint estimate (with $L^\infty_x$ replaced by $L^{2d/(d-2)}_x$ for scaling reasons) is now true, see \cite{tao:keel}.  It is also not hard to show that one can recover the estimate with only a logarithmic loss after compactifying time and frequency:

\begin{proposition}  Let $P$ be a Fourier multiplier with symbol in $C^\infty_0(\R^2)$, and let $I \subset \R$ be an interval.  Then
$$ \| e^{it\Delta} P f \|_{L^2_t L^\infty_x(I \times \R^2)} \lesssim_P \log(2+|I|)^{1/2} \|f\|_{L^2_x(\R^2)}$$
for all $f \in L^2_x(\R^2)$.
\end{proposition}

\begin{proof}(Sketch)  By the $TT^*$ method it suffices to show that
$$ |\int_I \int_I \langle P^* e^{i(t-t')\Delta} P F(t), F(t') \rangle_{L^2_x(\R^2)}\ dt dt'|
\lesssim_P \log(2+|I|) \|F\|_{L^2_t L^1_x(\R \times \R^2)}^2.$$
However, stationary phase computations yield the bound
$$ |\langle P^* e^{i(t-t')\Delta} P F(t), F(t') \rangle_{L^2_x(\R^2)}| \lesssim_P \frac{1}{1+|t-t'|}
\|F(t)\|_{L^1_x(\R^2)} \|F(t')\|_{L^1_x(\R^2)}.$$
The claim then follows from Schur's test or Young's inequality.
\end{proof}

We remark that one can achieve a similar result by removing the $P$ but instead placing $f$ in a Sobolev space $H^\eps_x(\R^2)$ for some $\eps > 0$; we omit the details.

In light of these near-misses, one might hope that some version of the endpoint Sobolev inequality without logarithmic losses could be salvaged.  One approach would be to pass from linear estimates to bilinear estimates, imposing some frequency separation on the two inputs; such estimates have proven to be rather useful in the study of nonlinear Schr\"odinger equations (as well as other for dispersive and wave models).  For instance, if $P, P'$ are Fourier multipliers with symbols in $C^\infty_0(\R^2)$ and with supports separated by a non-zero distance, then one has the estimate
$$ \| (e^{it\Delta} P f) (e^{it\Delta} P g) \|_{L^q_{t,x}(\R \times \R^2)} \lesssim_{P,P',q}
\|f\|_{L^2_x(\R^2)} \|g\|_{L^2_x(\R^2)}$$
for all $\frac{5}{3} < q \leq \infty$ and $f,g \in L^2_x(\R^2)$; see \cite{tao}.  This improves over what one can do just from \eqref{etf} and Bernstein's inequality, which can only handle the case $q \geq 2$.  This bilinear estimate is known to fail when $q < \frac{5}{3}$ but the endpoint $q = \frac{5}{3}$ remains open; see \cite{tao} for further discussion.  However, when $f,g,P,P'$ are spherically symmetric one can improve the range further; see \cite{shao}.

The main result of this paper is that while the bilinear setting undoubtedly improves the exponents in the non-endpoint case, it unfortunately does not do so in the endpoint case.

\begin{theorem}[Bilinear endpoint Strichartz fails]\label{main} Let $P, P'$ be Fourier multipliers with $C^\infty_0$ symbols which do not vanish identically.  Then there does not exist a constant $C > 0$ for which one has the estimate
$$ \| (e^{it\Delta} P f) (e^{it\Delta} P' g) \|_{L^1_t L^\infty_x(\R \times \R^2)} \leq C
\|f\|_{L^2_x(\R^2)} \|g\|_{L^2_x(\R^2)}$$
for all $f, g \in L^2_x(\R^2)$.
\end{theorem}

Note that this estimate would have followed from the linear endpoint Strichartz estimate by H\"older's inequality, if that estimate was true.

The proof of Theorem \ref{main} turns out to be remarkably ``low-tech'', and proceeds by using Theorem \ref{counter} as a ``black box''.  The basic idea is to remove the effect of $g$ by a standard randomised sign argument, thus reducing Theorem \ref{main} to Theorem \ref{counter}.

We thank Ioan Bejenaru for posing this question in the Schr\"odinger setting, and Sergiu Klainerman for posing it in the wave setting (see Section \ref{wave-sec} below).  The author is supported by a MacArthur Fellowship.

\section{Proof of Theorem \ref{main}}

Fix $P, P'$; we allow all implied constants in the $\lesssim$ notation to depend on these multipliers.  Suppose for contradiction that we did have an estimate
\begin{equation}\label{ty}
 \| (e^{it\Delta} P f) (e^{it\Delta} P' g) \|_{L^1_t L^\infty_x(\R \times \R^2)} \lesssim
\|f\|_{L^2_x(\R^2)} \|g\|_{L^2_x(\R^2)}
\end{equation}
for all $f, g \in L^2_x(\R^2)$.  For technical reasons it is convenient to exploit the frequency localisation (via the uncertainty principle) to replace the time axis $\R$ by the discretised variant $\Z$:

\begin{lemma}\label{disc}  We have
$$ \| (e^{it\Delta} P f) (e^{it\Delta} P' g) \|_{l^1_t L^\infty_x(\Z \times \R^2)} \lesssim
\|f\|_{L^2_x(\R^2)} \|g\|_{L^2_x(\R^2)}$$
for all $f, g \in L^2_x(\R^2)$.  
\end{lemma}

\begin{proof}  By the usual limiting arguments we may assume that $f,g$ are Schwartz functions (in order to justify all computations below).
We can find Fourier multipliers $\tilde P, \tilde P'$ with $C^\infty_0$ symbols such that $P = \tilde P P$ and $P' = \tilde P' P'$.  We then write
$$ (e^{it\Delta} P f) (e^{it\Delta} P' g) = \int_0^1 
[(e^{-i\theta \Delta} \tilde P) e^{i(t+\theta)\Delta} Pf]
[(e^{-i\theta \Delta} \tilde P') e^{i(t+\theta)\Delta} Pg]\ d\theta.$$
The convolution kernels $K_\theta(y)$, $K'_{\theta}(y)$ 
of $(e^{-i\theta \Delta} \tilde P)$ and $e^{-i\theta \Delta} \tilde P'$ are bounded
uniformly in magnitude by $\lesssim \frac{1}{1+|y|^{10}}$ (say).  Thus by Minkowski's inequality we have
\begin{align*} \| (e^{it\Delta} P f) (e^{it\Delta} P' g) \|_{L^\infty_x(\R^2)} &\lesssim
\int_0^1\int_{\R^2} \int_{\R^2}
\frac{1}{1+|y|^{10}} \frac{1}{1+|y'|^{10}}\\
\| (e^{i(t+\theta)\Delta} &P \tau_y f) (e^{i(t+\theta)\Delta} P \tau_{y'} g) \|_{L^\infty_x(\R^2)}\ dy dy' d\theta
\end{align*}
where $\tau_y$ is the operation of spatial translation by $y$.
Summing this in time and using Fubini's theorem one obtains
\begin{align*}\| (e^{it\Delta} P f) (e^{it\Delta} P' g) \|_{l^1_t L^\infty_x(\Z \times \R^2)} &\lesssim
\int_{\R^2} \int_{\R^2}
\frac{1}{1+|y|^{10}} \frac{1}{1+|y'|^{10}}\\
\| (e^{it\Delta} &P \tau_y f) (e^{it\Delta} P \tau_{y'} g) \|_{L^1_t L^\infty_x(\R \times \R^2)}\ dy dy'.
\end{align*}
Applying \eqref{ty} and using the integrability of $\frac{1}{1+|y|^{10}}$, the claim follows.
\end{proof}

Fix $N \geq 0$ and an arbitrary sequence of points $(x_n)_{n =-N}^N$ in $\R^2$.  From the above lemma we see that
$$ \sum_{n=-N}^N |e^{in\Delta} P f(x_n)| |e^{in\Delta} P' g(x_n)| \lesssim \|f\|_{L^2_x(\R^2)} \|g\|_{L^2_x(\R^2)}$$
for all Schwartz $f,g$.
Let us fix $f$ and dualise the above estimate in $g$, to obtain
$$ \| \sum_{n=-N}^N |e^{in\Delta} P f(x_n)| \epsilon_n (P')^* e^{-in\Delta} \delta_{x_n} \|_{L^2_x(\R^2)}
\lesssim \|f\|_{L^2_x(\R^2)} $$
for any $N \geq 0$, and any sequence $(\epsilon_n)_{n=-N}^N$ of scalars bounded in magnitude by $1$, where $(P')^*$ is the adjoint of $P'$ and $\delta_{x_n}$ is the Dirac mass at $x_n$.  We 
specialise $\epsilon_n = \pm 1$ to be iid random signs, square both sides, and take expectations (or use Khinchine's inequality) to obtain
$$ (\sum_{n=-N}^N \| |e^{in\Delta} P f(x_n)| (P')^* e^{-in\Delta} \delta_{x_n} \|_{L^2_x(\R^2)}^2)^{1/2}
\lesssim \|f\|_{L^2_x(\R^2)}.$$
But from Plancherel's theorem (and the hypothesis that the symbol of $P'$ does not vanish identically) we see that
$$ \| |e^{in\Delta} P f(x_n)| (P')^* e^{-in\Delta} \delta_{x_n} \|_{L^2_x(\R^2)} \gtrsim |e^{in\Delta} P f(x_n)| $$
so we conclude that
$$ (\sum_{n=-N}^N |e^{in\Delta} P f(x_n)|^2)^{1/2}
\lesssim \|f\|_{L^2_x(\R^2)}.$$
Since $x_n$ and $N$ were arbitrary, standard limiting arguments thus give us
$$ \| e^{in\Delta} Pf \|_{l^2_n L^\infty_x(\Z \times \R^2)} \lesssim \| f\|_{L^2_x(\R^2)}$$
for all Schwartz $f$.  Applying the unitary operator $e^{i\theta\Delta} f$ for $\theta \in [0,1]$ (which commutes with $e^{in\Delta} P$) and then averaging in $L^2_\theta$ then gives us
$$ \| e^{it\Delta} Pf \|_{L^2_t L^\infty_x(\R \times \R^2)} \lesssim \| f\|_{L^2_x(\R^2)}.$$
But this contradicts Theorem \ref{counter}, and we are done.
\endprf

\section{A variant for the wave equation}\label{wave-sec}

Observe that the above argument used very little about the Schr\"odinger propagators $e^{it\Delta}$, other than the group law and the fact that the convolution kernel of such propagators was uniformly rapidly decreasing  once one localised in both time and frequency.  One can thus adapt the above argument to other multipliers such as the wave propagators $e^{\pm i \sqrt{-\Delta}}$ in three dimensions.  The analogue of Theorem \ref{counter} is then

\begin{theorem}[Linear endpoint Strichartz fails]\label{counter-2}\cite{montgomery-smith}, \cite{tao-ex} Let $\epsilon = \pm 1$ be a sign, and let $P$ be a Fourier multiplier whose symbol lies in $C^\infty_0$, vanishes near the origin, and is not identically zero.
There does not exist a constant $C > 0$ for which one has the estimate
$$ \| e^{\epsilon it\sqrt{-\Delta}} P f \|_{L^2_t L^\infty_x(\R \times \R^3)} \leq C\|f\|_{L^2_x(\R^3)}$$
for all $f \in L^2_x(\R^3)$.
\end{theorem}

A routine modification of the above arguments now reveals that the corresponding bilinear estimate also fails:

\begin{corollary}[Biinear endpoint Strichartz fails] Let $\epsilon, \epsilon' = \pm 1$ be signs, and let $P$, $P'$ be Fourier multipliers whose symbol lies in $C^\infty_0$, vanishes near the origin, and is not identically zero.
There does not exist a constant $C > 0$ for which one has the estimate
$$ \| (e^{\epsilon it\sqrt{-\Delta}} P f) (e^{\epsilon' it\sqrt{-\Delta}} P' g) \|_{L^1_t L^\infty_x(\R \times \R^3)} \leq C\|f\|_{L^2_x(\R^2)} \|g\|_{L^2_x(\R^2)}$$
for all $f, g \in L^2_x(\R^2)$.
\end{corollary}

As there are no new ingredients in the proof we omit the details.

\section{Appendix: proof of Theorem \ref{counter}}

We now sketch the proof of Theorem \ref{counter}, following the Brownian motion ideas of \cite{montgomery-smith}.

Suppose for contradiction that Theorem \ref{counter} failed.  
Then by repeating the arguments used to prove Lemma \ref{disc} we have
$$ \| e^{it\Delta} P f \|_{l^2_t L^\infty_x(\Z \times \R^2)} \lesssim_P \|f\|_{L^2_x(\R^2)}.$$
Dualising as before, we are eventually obtain
$$ \| \sum_{n=-N}^N c_n P^* e^{-in\Delta} \delta_{x_n} \|_{L^2_x(\R^2)}
\lesssim (\sum_{n=-N}^N |c_n|^2)^{1/2}$$
for any $N \geq 0$, any points $(x_n)_{n=-N}^N$ in $\R^2$, and any complex numbers $(c_n)_{n=-N}^N$.  We set $c_n \equiv 1$ and then square to obtain
\begin{equation}\label{spin}
 |\sum_{n=-N}^N \sum_{n'=-N}^N  \langle P P^* e^{-i(n-n')\Delta} \delta_{x_n}, \delta_{x_{n'}} \rangle|
\lesssim N.
\end{equation}
Now we define $x_n$ by a random walk, so that $x_0 := 0$ and $x_{n+1}-x_n$ are i.i.d. Gaussian variables with variance $\sigma^2$ for some $\sigma > 0$ to be chosen later.  Standard probability theory then implies that $x_{n'}-x_n$ has a Gaussian distribution with variance $|n'-n| \sigma^2$.  By choosing $\sigma$ appropriately small but non-zero, and using the explicit formula for the convolution kernel of $e^{-i(n-n')\Delta}$ (and the fact that the convolution kernel of $PP^*$ is rapidly decreasing but has strictly positive integral) one can easily compute an expected lower bound
$$ \Re {\bf E} \langle P P^* e^{-i(n-n')\Delta} \delta_{x_n}, \delta_{x_{n'}} \rangle \gtrsim_{\sigma,P} \frac{1}{|n-n'|}$$
when $|n-n'|$ is larger than some constant $C_{\sigma,P}$ depending only on $\sigma$ and $P$.  Summing this (and using the crude bound of $O_{\sigma,P}(1)$ for the case when $n-n'$ is bounded) we obtain
$$ \Re {\bf E} \sum_{n=-N}^N \langle P P^* e^{-i(n-n')\Delta} \delta_{x_n}, \delta_{x_{n'}} \rangle \gtrsim_{\sigma,P} 
N \log N - O_{\sigma,P}(N)$$
which contradicts \eqref{spin} if $N$ is taken sufficiently large depending on $\sigma$ and $P$.
\endprf

\end{document}